\documentclass[11pt]{article}
\usepackage[russian]{babel}
\usepackage[cp1251]{inputenc}
\usepackage{amsmath,latexsym}
\usepackage[psamsfonts]{amssymb}
\usepackage[dvips]{graphicx}
\usepackage[mathcal]{euscript}
\begin{document}

MSC2010.35G31

\begin{center}
\textbf{Solution of the problem with initial conditions for the equation \\
high-order fractional Hilfer derivative
 }\\
\textbf{B.Yu.Irgashev}\\
\emph{Namangan  Engineering Construction Institute, \\ Namangan Branch of the Institute of Mathematics of the Republic of Uzbekistan.}\\
  \emph{E-mail: bahromirgasev@gmail.com}
\end{center}

\textbf{Abstract.} \emph{ In this paper, with the help of previously constructed self-similar solutions, we construct a solution to a Cauchy-type problem for an even-order high-order equation with a fractional derivative in the sense of Hilfer.}

 \textbf{Keywords.} \emph{ Higher order equation, Hilfer fractional derivative, self-similar solution, Cauchy problem.}

\begin{center}
\textbf{1. Introduction and construction of self-similar solutions
}\end{center}
Consider a degenerate fractional-order equation in the sense of Hilfer
\[{x^m}D_{0y}^{{\alpha _1},{\beta _1}}u\left( {x,y} \right) - d{y^k}D_{0x}^{{\alpha _2},{\beta _2}}u\left( {x,y} \right) = 0,\eqno(1)\]
here
\[m,k \in R,\,\,m >  - \,{\alpha _2};\,k >  - {\alpha _1}\,;\,q - 1 < {\alpha _1} \le q,\]
\[p - 1 < {\alpha _2} \le p,\,\,q < p,\,\,q,p \in N,\]
\[0 \le {\beta _1} \le 1,\,\,0 \le {\beta _2} \le 1,\,d =  \pm 1,\]
$D_{0t}^{\alpha ,\beta }$ fractional differentiation operator in the sense of Hilfer of order$\alpha$ and type $\beta$ (see [1])
\[D_{0t}^{\alpha ,\beta }u\left( {x,y} \right) = I_{oy}^{\beta \left( {s - \alpha } \right)}\frac{{{\partial ^s}}}{{\partial {t^s}}}\left( {I_{oy}^{\left( {1 - \beta } \right)\left( {s - \alpha } \right)}u\left( {x,y} \right)} \right),\]
\[s - 1 < \alpha  \le s,\,0 \le \beta  \le 1,\,s \in N.\]
In [2], using the methods of special operators, self-similar solutions of an equation with constant coefficients containing fractional derivatives in the sense of Riemann-Liouville were found. For an equation with an ordinary derivative of the form
\[\frac{{{\partial ^p}u\left( {x,y} \right)}}{{\partial {x^p}}} - \frac{{{\partial ^q}u\left( {x,y} \right)}}{{\partial {y^q}}} = 0,\,p < q,\]
self-similar solutions were constructed in [3]. \\
Note that various initial-boundary value problems for an equation with a fractional Hilfer derivative were studied in [4] - [6] and others.
In this paper, we construct self-similar solutions by a method that does not require knowledge of special operators.

We seek self-similar solutions to equation (1) in the form of the series
\[u\left( {x,y} \right) = {y^b}\sum\limits_{n = 0}^\infty  {{c_n}{{\left( {{x^a}{y^s}} \right)}^{n + \gamma }}}  = \sum\limits_{n = 0}^\infty  {{c_n}{x^{an + a\gamma }}{y^{sn + s\gamma  + b}}} ,\eqno(2)\]
here the parameters ${c_n},a,b,s,\gamma  \in R  $  are to be determined. \\
We first introduce the notation
\[{\overline {\left( a \right)} _s} = a\left( {a - 1} \right)...\left( {a - \left( {s - 1} \right)} \right),\,\,{\overline {\left( a \right)} _0} = 1,\,\,{\overline {\left( a \right)} _1} = a.\]
Substitute (2) into (1). Then we will have
\[\frac{{{\partial ^q}}}{{\partial {y^q}}}\left( {I_{oy}^{\left( {1 - {\beta _1}} \right)\left( {q - {\alpha _1}} \right)}u\left( {x,y} \right)} \right) = \sum\limits_{n = 0}^\infty  {\frac{{{c_n}{x^{an + a\gamma }}}}{{\Gamma \left( {\left( {1 - {\beta _1}} \right)\left( {q - {\alpha _1}} \right)} \right)}}\frac{{{\partial ^q}}}{{\partial {y^q}}}} \int\limits_0^y {\frac{{{z^{sn + s\gamma  + b}}}}{{{{\left( {y - z} \right)}^{1 - \left( {1 - {\beta _1}} \right)\left( {q - {\alpha _1}} \right)}}}}dz}  = \]
\[ = \sum\limits_{n = 0}^\infty  {\frac{{{c_n}\Gamma \left( {sn + s\gamma  + b + 1} \right){x^{an + a\gamma }}}}{{\Gamma \left( {sn + s\gamma  + b + 1 + \left( {1 - {\beta _1}} \right)\left( {q - {\alpha _1}} \right)} \right)}}\frac{{{\partial ^q}{y^{sn + s\gamma  + b + \left( {1 - {\beta _s}} \right)\left( {q - {\alpha _s}} \right)}}}}{{\partial {y^q}}}}  = \]
\[ = \sum\limits_{n = 0}^\infty  {\frac{{{c_n}\Gamma \left( {sn + s\gamma  + b + 1} \right){x^{an + a\gamma }}{y^{kn + k\gamma  + b + \left( {1 - {\beta _1}} \right)\left( {q - {\alpha _1}} \right) - q}}}}{{\Gamma \left( {sn + s\gamma  + b + 1 + \left( {1 - {\beta _1}} \right)\left( {q - {\alpha _1}} \right) - q} \right)}}} ,\]
further
\[I_{oy}^{{\beta _1}\left( {q - {\alpha _1}} \right)}\left( {\sum\limits_{n = 0}^\infty  {\frac{{{c_n}\Gamma \left( {sn + s\gamma  + b + 1} \right){x^{an + a\gamma }}{y^{sn + s\gamma  + b + \left( {1 - {\beta _1}} \right)\left( {q - {\alpha _1}} \right) - q}}}}{{\Gamma \left( {sn + s\gamma  + b + 1 + \left( {1 - {\beta _1}} \right)\left( {q - {\alpha _1}} \right) - q} \right)}}} } \right) = \]
\[ = \sum\limits_{n = 0}^\infty  {\frac{{{c_n}\Gamma \left( {sn + s\gamma  + b + 1} \right){x^{an + a\gamma }}}}{{\Gamma \left( {{\beta _1}\left( {q - {\alpha _1}} \right)} \right)\Gamma \left( {sn + s\gamma  + b + 1 + \left( {1 - {\beta _1}} \right)\left( {q - {\alpha _1}} \right) - q} \right)}}} \int\limits_0^y {\frac{{{z^{sn + s\gamma  + b + \left( {1 - {\beta _1}} \right)\left( {q - {\alpha _1}} \right) - q}}dz}}{{{{\left( {y - z} \right)}^{1 - {\beta _1}\left( {q - {\alpha _1}} \right)}}}}}  = \]
\[ = \sum\limits_{n = 0}^\infty  {\frac{{{c_n}\Gamma \left( {sn + s\gamma  + b + 1} \right){x^{an + a\gamma }}{y^{sn + s\gamma  + b - {\alpha _1}}}}}{{\Gamma \left( {sn + s\gamma  + b - {\alpha _1} + 1} \right)}}}.\eqno(3)\]
Similarly
\[\frac{{{\partial ^p}}}{{\partial {x^p}}}\left( {I_{oy}^{\left( {1 - {\beta _2}} \right)\left( {p - {\alpha _2}} \right)}u\left( {x,y} \right)} \right) = \sum\limits_{n = 0}^\infty  {\frac{{{c_n}{y^{sn + s\gamma  + b}}}}{{\Gamma \left( {\left( {1 - {\beta _2}} \right)\left( {p - {\alpha _2}} \right)} \right)}}\frac{{{\partial ^p}}}{{\partial {x^p}}}} \int\limits_0^x {\frac{{{z^{an + a\gamma }}}}{{{{\left( {x - z} \right)}^{1 - \left( {1 - {\beta _2}} \right)\left( {p - {\alpha _2}} \right)}}}}dz}  = \]
\[ = \sum\limits_{n = 0}^\infty  {\frac{{{c_n}{y^{sn + s\gamma  + b}}}}{{\Gamma \left( {\left( {1 - {\beta _2}} \right)\left( {p - {\alpha _2}} \right)} \right)}}\frac{{{\partial ^p}{x^{an + a\gamma  + \left( {1 - {\beta _2}} \right)\left( {p - {\alpha _2}} \right)}}}}{{\partial {x^p}}}} \int\limits_0^1 {\frac{{{z^{an + a\gamma }}}}{{{{\left( {1 - z} \right)}^{1 - \left( {1 - {\beta _2}} \right)\left( {p - {\alpha _2}} \right)}}}}dz}  = \]
\[ = \sum\limits_{n = 0}^\infty  {\frac{{{{\overline {\left( {an + a\gamma  + \left( {1 - {\beta _2}} \right)\left( {p - {\alpha _2}} \right)} \right)} }_p}\Gamma \left( {an + a\gamma  + 1} \right){c_n}{x^{an + a\gamma  + \left( {1 - {\beta _2}} \right)\left( {p - {\alpha _2}} \right) - p}}{y^{sn + s\gamma  + b}}}}{{\Gamma \left( {an + a\gamma  + 1 + \left( {1 - {\beta _2}} \right)\left( {p - {\alpha _2}} \right)} \right)}}}  = \]
\[ = \sum\limits_{n = 0}^\infty  {\frac{{\Gamma \left( {an + a\gamma  + 1} \right){c_n}{x^{an + a\gamma  + \left( {1 - {\beta _2}} \right)\left( {p - {\alpha _2}} \right) - p}}{y^{sn + s\gamma  + b}}}}{{\Gamma \left( {an + a\gamma  + 1 + \left( {1 - {\beta _2}} \right)\left( {p - {\alpha _2}} \right) - p} \right)}}} ,\]
further
\[I_{ox}^{{\beta _2}\left( {p - {\alpha _2}} \right)}\left( {\sum\limits_{n = 0}^\infty  {\frac{{\Gamma \left( {an + a\gamma  + 1} \right){c_n}{x^{an + a\gamma  + \left( {1 - {\beta _2}} \right)\left( {p - {\alpha _2}} \right) - p}}{y^{sn + s\gamma  + b}}}}{{\Gamma \left( {an + a\gamma  + 1 + \left( {1 - {\beta _2}} \right)\left( {p - {\alpha _2}} \right) - p} \right)}}} } \right) = \]
\[ = \sum\limits_{n = 0}^\infty  {\frac{{\Gamma \left( {an + a\gamma  + 1} \right){c_n}{y^{sn + s\gamma  + b}}}}{{\Gamma \left( {an + a\gamma  + 1 + \left( {1 - {\beta _2}} \right)\left( {p - {\alpha _2}} \right) - p} \right)\Gamma \left( {{\beta _2}\left( {p - {\alpha _2}} \right)} \right)}}\int\limits_0^x {\frac{{{z^{an + a\gamma  + \left( {1 - {\beta _2}} \right)\left( {p - {\alpha _2}} \right) - p}}dz}}{{{{\left( {x - z} \right)}^{1 - {\beta _2}\left( {p - {\alpha _2}} \right)}}}}} }  = \]
\[ = \sum\limits_{n = 0}^\infty  {\frac{{\Gamma \left( {an + a\gamma  + 1} \right){c_n}{y^{sn + s\gamma  + b}}{x^{an + a\gamma  - {\alpha _2}}}}}{{\Gamma \left( {an + a\gamma  + 1 + \left( {1 - {\beta _2}} \right)\left( {p - {\alpha _2}} \right) - p} \right)\Gamma \left( {{\beta _2}\left( {p - {\alpha _2}} \right)} \right)}}\int\limits_0^1 {\frac{{{z^{an + a\gamma  + \left( {1 - {\beta _2}} \right)\left( {p - {\alpha _2}} \right) - p}}dz}}{{{{\left( {1 - z} \right)}^{1 - {\beta _2}\left( {p - {\alpha _2}} \right)}}}}} }  = \]
\[ = \sum\limits_{n = 0}^\infty  {\frac{{\Gamma \left( {an + a\gamma  + 1} \right){c_n}{y^{sn + s\gamma  + b}}{x^{an + a\gamma  - {\alpha _2}}}}}{{\Gamma \left( {an + a\gamma  - {\alpha _2} + 1} \right)}}}.\eqno(4)\]
Substituting (3), (4) into (1), we obtain
\[\sum\limits_{n = 0}^\infty  {\frac{{{c_n}\Gamma \left( {sn + s\gamma  + b + 1} \right){{\left( {{x^a}{y^s}} \right)}^{n + \gamma }}{x^m}{y^{ - {\alpha _1}}}}}{{\Gamma \left( {sn + s\gamma  + b - {\alpha _1} + 1} \right)}}}  = d\sum\limits_{n = 0}^\infty  {\frac{{\Gamma \left( {an + a\gamma  + 1} \right){c_n}{{\left( {{x^a}{y^s}} \right)}^{n + \gamma }}{y^k}{x^{ - {\alpha _2}}}}}{{\Gamma \left( {an + a\gamma  - {\alpha _2} + 1} \right)}}}  \Rightarrow \]
\[\sum\limits_{n = 0}^\infty  {\frac{{{c_n}\Gamma \left( {sn + s\gamma  + b + 1} \right){{\left( {{x^a}{y^s}} \right)}^{n + \gamma }}{x^{m + {\alpha _2}}}{y^{ - {\alpha _1} - k}}}}{{\Gamma \left( {sn + s\gamma  + b - {\alpha _1} + 1} \right)}}}  = d\sum\limits_{n = 0}^\infty  {\frac{{\Gamma \left( {an + a\gamma  + 1} \right){c_n}{{\left( {{x^a}{y^s}} \right)}^{n + \gamma }}}}{{\Gamma \left( {an + a\gamma  - {\alpha _2} + 1} \right)}}}  \Rightarrow \]
\[a = m + {\alpha _2},\,s =  - {\alpha _1} - k \Rightarrow \]
\[\sum\limits_{n = 0}^\infty  {\frac{{{c_n}\Gamma \left( {sn + s\gamma  + b + 1} \right){{\left( {{x^a}{y^s}} \right)}^{n + \gamma  + 1}}}}{{\Gamma \left( {sn + s\gamma  + b - {\alpha _1} + 1} \right)}}}  = d\sum\limits_{n = 0}^\infty  {\frac{{\Gamma \left( {an + a\gamma  + 1} \right){c_n}{{\left( {{x^a}{y^s}} \right)}^{n + \gamma }}}}{{\Gamma \left( {an + a\gamma  - {\alpha _2} + 1} \right)}}}  \Rightarrow \]
\[\sum\limits_{n = 0}^\infty  {{c_n}\frac{{\Gamma \left( {sn + s\gamma  + b + 1} \right){{\left( {{x^a}{y^s}} \right)}^{n + \gamma  + 1}}}}{{\Gamma \left( {sn + s\gamma  + b - {\alpha _1} + 1} \right)}}}  = d\sum\limits_{n = 0}^\infty  {{c_n}\frac{{{{\overline {\left( {an + a\gamma  - {\alpha _2} + p} \right)} }_p}\Gamma \left( {an + a\gamma  + 1} \right){{\left( {{x^a}{y^s}} \right)}^{n + \gamma }}}}{{\Gamma \left( {an + a\gamma  - {\alpha _2} + p + 1} \right)}}}  \Rightarrow \]
got conditions for the parameter $\gamma :$
\[{\overline {\left( {a\gamma  - {\alpha _2} + p} \right)} _p} = 0 \Rightarrow {\gamma _j} = \frac{{{\alpha _2} - j}}{a} = \frac{{{\alpha _2} - j}}{{{\alpha _2} + m}},\,\,j = 1,2,...,p.\]
Substituting the found parameters into representation (2), we find the formula for the coefficients ${c_n}$. We have
\[d\frac{{\Gamma \left( {an + a\gamma  + 1} \right)}}{{\Gamma \left( {an + a\gamma  - {\alpha _2} + 1} \right)}}{c_n} = {c_{n - 1}}\frac{{\Gamma \left( {s\left( {n - 1} \right) + s\gamma  + b + 1} \right)}}{{\Gamma \left( {s\left( {n - 1} \right) + s\gamma  + b - {\alpha _1} + 1} \right)}} \Rightarrow \]
\[{c_n} = {c_{n - 1}}d\frac{{\Gamma \left( {s\left( {n - 1} \right) + s\gamma  + b + 1} \right)\Gamma \left( {an + a\gamma  - {\alpha _2} + 1} \right)}}{{\Gamma \left( {s\left( {n - 1} \right) + s\gamma  + b - {\alpha _1} + 1} \right)\Gamma \left( {an + a\gamma  + 1} \right)}} \Rightarrow \]
\[{c_n} = \frac{{\Gamma \left( { - \left( {{\alpha _1} + k} \right)\left( {n - 1 + \frac{{\left( {{\alpha _2} - j} \right)}}{{{\alpha _2} + m}}} \right) + b + 1} \right)\Gamma \left( {\left( {m + {\alpha _2}} \right)n - j + 1} \right){c_{n - 1}}d}}{{\Gamma \left( { - \left( {{\alpha _1} + k} \right)\left( {n - 1 + \frac{{\left( {{\alpha _2} - j} \right)}}{{{\alpha _2} + m}}} \right) + b - {\alpha _1} + 1} \right)\Gamma \left( {\left( {m + {\alpha _2}} \right)n + {\alpha _2} - j + 1} \right)}} \Rightarrow \]
\[c_n^j = {c_0}{d^n}\prod\limits_{l = 1}^n {\frac{{\Gamma \left( { - \left( {{\alpha _1} + k} \right)\left( {l - 1 + \frac{{\left( {{\alpha _2} - j} \right)}}{{{\alpha _2} + m}}} \right) + b + 1} \right)\Gamma \left( {\left( {m + {\alpha _2}} \right)l - j + 1} \right)}}{{\Gamma \left( { - \left( {{\alpha _1} + k} \right)\left( {l - 1 + \frac{{\left( {{\alpha _2} - j} \right)}}{{{\alpha _2} + m}}} \right) + b - {\alpha _1} + 1} \right)\Gamma \left( {\left( {m + {\alpha _2}} \right)l + {\alpha _2} - j + 1} \right)}}}.\eqno(5)\]
So, self-similar solutions of equation (1) have the form
\[{u_j}\left( {x,y} \right) = {y^b}{t^{1 - \frac{{m + j}}{{{\alpha _2} + m}}}}\sum\limits_{n = 0}^\infty  {c_n^j{t^n}},\eqno(6)\]
where
\[t = {x^{m + {\alpha _2}}}{y^{ - {\alpha _1} - k}},j = 1,2,...,p.\]
Let's consider some special cases:

1). Let $ m = s = 0 $. Then from formula (5) we have
\[c_n^j = {c_0}{d^n}\prod\limits_{l = 1}^n {\frac{{\Gamma \left( { - {\alpha _1}\left( {l - \frac{j}{{{\alpha _2}}}} \right) + b + 1} \right)\Gamma \left( {{\alpha _2}l - j + 1} \right)}}{{\Gamma \left( { - {\alpha _1}\left( {l - \frac{j}{{{\alpha _2}}} + 1} \right) + b + 1} \right)\Gamma \left( {{\alpha _2}\left( {l + 1} \right) - j + 1} \right)}}}  = \]
\[c_n^j = {c_0}{d^n}\frac{{\Gamma \left( { - {\alpha _1}\left( {1 - \frac{j}{{{\alpha _2}}}} \right) + b + 1} \right)\Gamma \left( {{\alpha _2} - j + 1} \right)}}{{\Gamma \left( { - {\alpha _1}\left( {n - \frac{j}{{{\alpha _2}}} + 1} \right) + b + 1} \right)\Gamma \left( {{\alpha _2}\left( {n + 1} \right) - j + 1} \right)}},\]
if now
\[{c_0} = \frac{1}{{\Gamma \left( { - {\alpha _1}\left( {1 - \frac{j}{{{\alpha _2}}}} \right) + b + 1} \right)\Gamma \left( {{\alpha _2} - j + 1} \right)}},\]
then formula (6) will be represented as
\[{u_j}\left( {x,y} \right) = {y^b}{t^{1 - \frac{j}{{{\alpha _2}}}}}\sum\limits_{n = 0}^\infty  {\frac{{{{\left( {dt} \right)}^n}}}{{\Gamma \left( { - {\alpha _1}n + \frac{{{\alpha _1}j}}{{{\alpha _2}}} - {\alpha _1} + b + 1} \right)\Gamma \left( {{\alpha _2}n + {\alpha _2} - j + 1} \right)}}}=\]
\[ = {y^b}{t^{1 - \frac{j}{{{\alpha _2}}}}}{W_{\left( { - {\alpha _1}, - {\alpha _1} + \frac{{{\alpha _1}}}{{{\alpha _2}}}j + b + 1} \right),\left( {{\alpha _2},{\alpha _2} - j + 1} \right)}}\left( {dt} \right),\eqno(7)\]
here
\[{W_{\left( {\mu ,a} \right),\left( {\nu ,b} \right)}}\left( z \right) = \sum\limits_{n = 0}^\infty  {\frac{{{z^n}}}{{\Gamma \left( {\mu n + a} \right)\Gamma \left( {\nu n + b} \right)}}} ,\,\mu ,\nu  \in R,\,a,b \in C,\,\mu  + \nu  > 0\]
- generalized Wright function [2].
Series (7) converges uniformly and it can be differentiated term by term, since $ {\alpha_2}> {\alpha _1} $.
Note that there is no condition on the $ b $ parameter. Note that solutions of the form (7) exactly coincide with the solutions obtained in [2] for the equation with the fractional Riemann-Liouville derivative.

2). Now let in representation (7): $ {\alpha _2} = p \in N, \, \, {\alpha _1} = \alpha, \, \, q - 1 <\alpha \le q \in N, \, \, q <p $, then solutions of the equation
\[D_{0y}^{\alpha ,\beta }u\left( {x,y} \right) - d\frac{{{\partial ^p}u\left( {x,y} \right)}}{{\partial {x^p}}} = 0,\,\,\,0 \le \beta  \le 1,\eqno(8)\]
there will be expressions like
\[{u_s}\left( {x,y} \right) = {y^b}\sum\limits_{n = 0}^\infty  {\frac{{{d^n}{{\left( {x{y^{ - \frac{\alpha }{p}}}} \right)}^{pn + s}}}}{{\Gamma \left( { - \alpha n - \frac{\alpha }{p}s + b + 1} \right)\left( {pn + s} \right)!}}} ,\,\,s = 0,1,...,p - 1,\]
hence, their linear combination is also a solution to equation (8)
\[u\left( {x,y} \right) = {y^b}\left( {{c_0}\sum\limits_{n = 0}^\infty  {\frac{{{d^n}{t^{pn}}}}{{\left( {pn} \right)!\Gamma \left( { - \alpha n + b + 1} \right)}}}  + {c_1}\sum\limits_{n = 0}^\infty  {\frac{{{d^n}{t^{pn + 1}}}}{{\left( {pn + 1} \right)!\Gamma \left( { - \alpha \left( {n + \frac{1}{p}} \right) + b + 1} \right)}} + ...} } \right.\]
\[ + \left. {{c_{p - 1}}\sum\limits_{n = 0}^\infty  {\frac{{{d^n}{t^{pn + p - 1}}}}{{\left( {p\left( {n + 1} \right) - 1} \right)!\Gamma \left( { - \alpha \left( {n + \frac{{p - 1}}{p}} \right) + b + 1} \right)}}} } \right),\]
where $ {c_i} - $ arbitrary real numbers, $i=0,1,...,p-1 .$\\
Now let $ {c_i} $ be such that $ {c_i} \ne {c_j} $, for $ i \ne j $ and $ c_i ^ p = d $. Then we get
\[u\left( {x,y} \right) = {y^b}\phi \left( { - \frac{\alpha }{p},b + 1,cx{y^{ - \frac{\alpha }{p}}}} \right),\,{c^p} = d,\eqno(9)\]
here
\[\phi \left( { - \delta ,\varepsilon ,z} \right) = \sum\limits_{k = 0}^\infty  {\frac{{{z^k}}}{{k!\Gamma \left( { - \delta k + \varepsilon } \right)}}}\]
- Wright function. \\
\begin{center}
\textbf{2. Problem with initial conditions
}\end{center}
In this section, using the found particular solutions (9), we obtain an explicit form of the solution to the Cauchy problem for an equation of the form:
\[D_{0y}^{\alpha ,\beta }u\left( {x,y} \right) - {\left( { - 1} \right)^{n - 1}}\frac{{{\partial ^{2n}}u\left( {x,y} \right)}}{{\partial {x^{2n}}}} = 0,\,n=2,3,...,\,\,0 < \alpha  < 2.\eqno(10)\]
Note that the Cauchy problem for higher-order equations with Dzhrbashyan-Nersesyan and Riemann-Liouville fractional derivatives, orders $ 0 <\alpha <1 $, were considered, respectively, in [7], [8]. We will use ideas from these works.

Following (10), we have
\[c_{1k}^{2n} = {\left( { - 1} \right)^{n - 1}} \Rightarrow {c_{1k}} = {e^{\frac{{n - 1 - 2k}}{{2n}}i\pi }},\,k = \overline {0,\left( {n - 1} \right)} ,\,{\mathop{\rm Re}\nolimits} \,{c_{1k}} > 0,\]
Consider the function
\[{\Gamma _b}\left( {x - \xi ,y - \eta } \right) = \left\{ \begin{array}{l}
\Gamma _b^1\left( {x - \xi ,y - \eta } \right),\,x > \xi ,\\
\Gamma _b^2\left( {x - \xi ,y - \eta } \right),\,x < \xi ,
\end{array} \right.\eqno(11)\]
where
\[\Gamma _b^1\left( {x - \xi ,y - \eta } \right) = \frac{{{{\left( {y - \eta } \right)}^b}}}{{2n}}\sum\limits_{k = 0}^{n - 1} {\left( { - {e^{\frac{{n - 1 - 2k}}{{2n}}i\pi }}} \right)\phi \left( { - \frac{\alpha }{{2n}},b + 1, - {e^{\frac{{n - 1 - 2k}}{{2n}}i\pi }}t} \right),\,} \]
\[\Gamma _b^2\left( {x - \xi ,y - \eta } \right) =  - \frac{{{{\left( {y - \eta } \right)}^b}}}{{2n}}\sum\limits_{k = 0}^{n - 1} {{e^{\frac{{n - 1 - 2k}}{{2n}}i\pi }}\phi \left( { - \frac{\alpha }{{2n}},b + 1, - {e^{\frac{{n - 1 - 2k}}{{2n}}i\pi }}\left( { - t} \right)} \right)}, \]
\[t = \frac{{\left( {x - \xi } \right)}}{{{{\left( {y - \eta } \right)}^{\frac{\alpha }{{2n}}}}}} .\]
The following lemma is true (see [8]).

\textbf{Lemma 1.} For $ s \in N $, the following relation holds:
\[{\left. {\frac{{{\partial ^s}\Gamma _b^1\left( {x - \xi ,y - \eta } \right)}}{{\partial {x^s}}}} \right|_{x = \xi }} - {\left. {\frac{{{\partial ^s}\Gamma _b^2\left( {x - \xi ,y - \eta } \right)}}{{\partial {x^s}}}} \right|_{x = \xi }} = \]
\[ = {\left( { - 1} \right)^{n - 1}}\frac{{{{\left( {y - \eta } \right)}^{b - \frac{\alpha }{{2n}}s}}}}{{\Gamma \left( {b + 1 - \frac{\alpha }{{2n}}} \right)}}\left\{ \begin{array}{l}
1,\,\,s = \left( {2n - 1} \right)\,\left( {\bmod \,2n} \right),\\
0,\,\,s \ne \left( {2n - 1} \right)\,\left( {\bmod \,2n} \right).
\end{array} \right.\]
In what follows, we need the asymptotic behavior of the Wright function for large values of the variable. The main results on asymptotics were obtained by Wright [see 9]. In particular, the following theorem is proved.

\textbf{Theorem 1} [см.3].If $\left| {\arg \,y} \right| \le \min \left\{ {\frac{3}{2}\pi \left( {1 - \sigma } \right),\pi } \right\} - \varepsilon ,\,\,0 < \sigma  < 1,$ then
\[\phi \left( { - \sigma ,\beta ,t} \right) = {Y^{\frac{1}{2} - \beta }}{e^{ - Y}}\left\{ {\sum\limits_{m = 0}^{M - 1} {{A_m}{Y^{ - m}}}  + O\left( {{Y^{ - M}}} \right)} \right\},\,\,\left| t \right| \to \infty ,\eqno(12)\]
here $ Y = \left( {1 - \sigma } \right){\left( {{\sigma ^\sigma }y} \right)^{\frac{1}{{1 - \sigma }}}},\,\,y =  - t,\,$ coefficients $ {A_m} $ depends on $\sigma ,\beta $.\\
If $ n = 3,4, ... $, then $ 0 <\frac {\alpha} {{2n}} <\frac {1} {n} \le \frac {1} {3}, $ and the relation holds\[\left|{\arg \,y} \right| = \left| {\arg \,\left( { - t} \right)} \right| = \left| {\frac{{n - 1 - 2k}}{{2n}}\pi } \right| \le \pi  - \varepsilon , k=0,1,...,n-1.\eqno(13)\]
If $ n = 2 $, then for $ 0 <\frac {\alpha} {4} \le \frac {1} {3}, $ we have\[\left| {\arg \,y} \right| = \left| {\arg \,\left( { - t} \right)} \right| = \left| { \pm \frac{\pi }{4}} \right| < \pi ,\eqno(14)\]
and for $ \frac {1} {3} <\frac {\alpha} {4} <\frac {1} {2}, $ we get
\[\left| {\arg \,y} \right| = \left| {\arg \,\left( { - t} \right)} \right| = \left| { \pm \frac{\pi }{4}} \right| < \frac{3}{2}\pi \left( {1 - \frac{\alpha }{4}} \right).\eqno(15)\]
Taking into account (13) - (15), we conclude that for $ n = 2,3, ... $ we always have relation (12). Now we write (11) in the form
\[{\Gamma _b}\left( {x - \xi ,y - \eta } \right) = \frac{{{{\left( {y - \eta } \right)}^b}}}{{2n}}\sum\limits_{k = 0}^{n - 1} {\left( { - {e^{\frac{{n - 1 - 2k}}{{2n}}i\pi }}} \right)\phi \left( { - \frac{\alpha }{{2n}},b + 1, - {e^{\frac{{n - 1 - 2k}}{{2n}}i\pi }}\left| t \right|} \right),\,}\eqno(16) \]
where
\[\left| t \right| = \frac{{\left| {x - \xi } \right|}}{{{{\left( {y - \eta } \right)}^{\frac{\alpha }{{2n}}}}}},\]
then, taking into account (12), (13) - (15), for large values of $ \left | t \right | $, for (15) the estimate obtained in [8] is valid:
\[\left| {\phi \left( { - \frac{\alpha }{{2n}},b + 1, - {e^{\frac{{n - 1 - 2k}}{{2n}}i\pi }}\left| t \right|} \right)} \right| \le C{\left| t \right|^{\frac{{ - 2n}}{{2n - \alpha }}\left( {b + \frac{1}{2}} \right)}}\exp \left( { - \sigma {{\left| t \right|}^{\frac{{2n}}{{2n - \alpha }}}}} \right),\eqno(17)\]
here
\[\sigma  = \left( {1 - \frac{\alpha }{{2n}}} \right){\left( {\frac{\alpha }{{2n}}} \right)^{\frac{\alpha }{{2n - \alpha }}}}\cos \frac{{n - 1}}{{2n - \alpha }}\pi ,\,0 < C - const.\]
Let's find the following expression
\[\int\limits_{ - \infty }^{ + \infty } {\frac{{{\partial ^k}}}{{\partial {y^k}}}\left( {I_{0y}^{\left( {1 - \beta } \right)\left( {s - \alpha } \right)}{\Gamma _b}\left( {x - \xi ,y} \right)} \right)d\xi } ,\]
where
\[k = \overline {0,\left( {s - 1} \right)} ,\,s = 1,2.\]
We have
\[\frac{{{\partial ^k}}}{{\partial {y^k}}}\left( {I_{0y}^{\left( {1 - \beta } \right)\left( {s - \alpha } \right)}{\Gamma _b}\left( {x - \xi ,y} \right)} \right) = \]
\[ = \frac{1}{{\Gamma \left( {\left( {1 - \beta } \right)\left( {s - \alpha } \right)} \right)}}\frac{{{\partial ^k}}}{{\partial {y^k}}}\int\limits_0^y {\frac{{\sum\limits_{l = 0}^{n - 1} {\frac{{{z^b}}}{{2n}}\left( { - {\lambda _l}} \right)\phi \left( { - \frac{\alpha }{{2n}},b + 1, - {\lambda _l}\left| {x - \xi } \right|{z^{ - \frac{\alpha }{{2n}}}}} \right)} }}{{{{\left( {y - z} \right)}^{1 - \left( {1 - \beta } \right)\left( {s - \alpha } \right)}}}}dz = } \]
\[ = \sum\limits_{l = 0}^{n - 1} {\sum\limits_{m = 0}^\infty  {\frac{{\left( { - {\lambda _l}} \right){{\left( { - {\lambda _l}\left| {x - \xi } \right|} \right)}^m}{y^{ - \frac{\alpha }{{2n}}m + b + \left( {1 - \beta } \right)\left( {s - \alpha } \right) - k}}}}{{2nm!\Gamma \left( { - \frac{\alpha }{{2n}}m + b + \left( {1 - \beta } \right)\left( {s - \alpha } \right) - k + 1} \right)}}} } ,\]
\[\int\limits_{ - \infty }^{ + \infty } {\sum\limits_{l = 0}^{n - 1} {\sum\limits_{m = 0}^\infty  {\frac{{\left( { - {\lambda _l}} \right){{\left( { - {\lambda _l}\left| {x - \xi } \right|} \right)}^m}{y^{ - \frac{\alpha }{{2n}}m + b + \left( {1 - \beta } \right)\left( {s - \alpha } \right) - k}}}}{{2nm!\Gamma \left( { - \frac{\alpha }{{2n}}m + b + \left( {1 - \beta } \right)\left( {s - \alpha } \right) - k + 1} \right)}}} } d\xi }  = \]
\[ = \int\limits_{ - \infty }^x {\sum\limits_{l = 0}^{n - 1} {\sum\limits_{m = 0}^\infty  {\frac{{\left( { - {\lambda _l}} \right){{\left( { - {\lambda _l}\left( {x - \xi } \right)} \right)}^m}{y^{ - \frac{\alpha }{{2n}}m + b + \left( {1 - \beta } \right)\left( {s - \alpha } \right) - k}}}}{{2nm!\Gamma \left( { - \frac{\alpha }{{2n}}m + b + \left( {1 - \beta } \right)\left( {s - \alpha } \right) - k + 1} \right)}}} } d\xi }  + \]
\[ + \int\limits_x^{ + \infty } {\sum\limits_{l = 0}^{n - 1} {\sum\limits_{m = 0}^\infty  {\frac{{\left( { - {\lambda _l}} \right){{\left( { - {\lambda _l}\left( {\xi  - x} \right)} \right)}^m}{y^{ - \frac{\alpha }{{2n}}m + b + \left( {1 - \beta } \right)\left( {s - \alpha } \right) - k}}}}{{2nm!\Gamma \left( { - \frac{\alpha }{{2n}}m + b + \left( {1 - \beta } \right)\left( {s - \alpha } \right) - k + 1} \right)}}} } d\xi }  = \]
\[ = {J_1} + {J_2},\]
we calculate each term separately, taking into account estimate (17), we have
\[{J_1} = \int\limits_{ - \infty }^x {\sum\limits_{l = 0}^{n - 1} {\sum\limits_{m = 0}^\infty  {\frac{{\left( { - {\lambda _l}} \right){{\left( { - {\lambda _l}\left( {x - \xi } \right)} \right)}^m}{y^{ - \frac{\alpha }{{2n}}m + b + \left( {1 - \beta } \right)\left( {s - \alpha } \right) - k}}}}{{2nm!\Gamma \left( { - \frac{\alpha }{{2n}}m + b + \left( {1 - \beta } \right)\left( {s - \alpha } \right) - k + 1} \right)}}} } d\xi }  = \]
\[\left. {\sum\limits_{l = 0}^{n - 1} {\sum\limits_{m = 0}^\infty  {\frac{{{{\left( { - {\lambda _l}} \right)}^{m + 1}}{y^{ - \frac{\alpha }{{2n}}m + b + \left( {1 - \beta } \right)\left( {s - \alpha } \right) - k}}{{\left( {x - \xi } \right)}^{m + 1}}}}{{2n\left( {m + 1} \right)!\Gamma \left( { - \frac{\alpha }{{2n}}m + b + \left( {1 - \beta } \right)\left( {s - \alpha } \right) - k + 1} \right)}}} } } \right|_{\xi  =  - \infty }^{\xi  = x} = \]
\[ = \frac{{{y^{\frac{\alpha }{{2n}} + b + \left( {1 - \beta } \right)\left( {s - \alpha } \right) - k}}}}{{2n}}\sum\limits_{l = 0}^{n - 1} {\left( {\phi \left( { - \frac{\alpha }{{2n}},\frac{\alpha }{{2n}} + b + \left( {1 - \beta } \right)\left( {s - \alpha } \right) - k + 1,\left( { - {\lambda _l}x{y^{ - \frac{\alpha }{{2n}}}}} \right)} \right)} \right.}  - \]
\[\left. {\left. { - \frac{1}{{\Gamma \left( {\frac{\alpha }{{2n}} + b + \left( {1 - \beta } \right)\left( {s - \alpha } \right) - k + 1} \right)}}} \right)} \right|_{\xi  =  - \infty }^{\xi  = x} = \]
\[ = \frac{{{y^{\frac{\alpha }{{2n}} + b + \left( {1 - \beta } \right)\left( {s - \alpha } \right) - k}}}}{{2\Gamma \left( {\frac{\alpha }{{2n}} + b + \left( {1 - \beta } \right)\left( {s - \alpha } \right) - k + 1} \right)}}.\]
Similarly
\[{J_2} = \int\limits_x^{ + \infty } {\sum\limits_{l = 0}^{n - 1} {\sum\limits_{m = 0}^\infty  {\frac{{\left( { - {\lambda _l}} \right){{\left( { - {\lambda _l}\left( {\xi  - x} \right)} \right)}^m}{y^{ - \frac{\alpha }{{2n}}m + b + \left( {1 - \beta } \right)\left( {s - \alpha } \right) - k}}}}{{2nm!\Gamma \left( { - \frac{\alpha }{{2n}}m + b + \left( {1 - \beta } \right)\left( {s - \alpha } \right) - k + 1} \right)}}} } d\xi }  = \]
\[ = \frac{{{y^{\frac{\alpha }{{2n}} + b + \left( {1 - \beta } \right)\left( {s - \alpha } \right) - k}}}}{{2\Gamma \left( {\frac{\alpha }{{2n}} + b + \left( {1 - \beta } \right)\left( {s - \alpha } \right) - k + 1} \right)}}.\]
So
\[\int\limits_{ - \infty }^{ + \infty } {\frac{{{\partial ^k}}}{{\partial {y^k}}}\left( {I_{0y}^{\left( {1 - \beta } \right)\left( {s - \alpha } \right)}{\Gamma _b}\left( {x - \xi ,y} \right)} \right)d\xi }  = \frac{{{y^{\frac{\alpha }{{2n}} + b + \left( {1 - \beta } \right)\left( {s - \alpha } \right) - k}}}}{{\Gamma \left( {\frac{\alpha }{{2n}} + b + \left( {1 - \beta } \right)\left( {s - \alpha } \right) - k + 1} \right)}}.\eqno(18)\]
Taking into account (18), we can show the correctness of the following lemma.

\textbf{Lemma 2.} For any function $ h \left (x \right) \in C \left (R \right) $ such that
\[\left| {h\left( x \right)} \right| \le M\exp \left( {c{{\left| x \right|}^{\frac{{2n}}{{2n - \alpha }}}}} \right),\,\,c < \sigma ,\,\,\left| x \right| \to \infty,0<M-const,\]
equality holds
\[\mathop {\lim }\limits_{y \to  + 0} \int\limits_{ - \infty }^{ + \infty } {h\left( \xi  \right)\frac{{{\partial ^k}}}{{\partial {y^k}}}\left( {I_{0y}^{\left( {1 - \beta } \right)\left( {s - \alpha } \right)}{\Gamma _{ - \frac{\alpha }{{2n}} - \left( {1 - \beta } \right)\left( {s - \alpha } \right) + j}}\left( {x - \xi ,y} \right)} \right)d\xi }  = \]
\[ = h\left( x \right) \cdot \left\{ \begin{array}{l}
0,\,j \ne k,\\
1,\,\,j = k,
\end{array} \right.\]
here
\[s = 1,2;\,s - 1 < \alpha  \le s;\,k = \overline {0,s - 1} ;\,0 \le \beta  \le 1;\,j \in N.\]

\textbf{Cauchy problem.} Find the solution $ u \left ({x, y} \right) $ to equation (10)
in the region $ D = \left \{{\left ({x, y} \right): - \infty <x <+ \infty, 0 <y} \right \}, $ satisfying the following conditions:

1) $D_{0y}^{\alpha ,\beta }u\left( {x,y} \right),\,\frac{{{\partial ^{2n}}u\left( {x,y} \right)}}{{\partial {x^{2n}}}} \in C\left( D \right);$

2) $I_{0y}^{\left( {1 - \beta } \right)\left( {s - \alpha } \right)}u\left( {x,y} \right) \in {C^k}\left( {\overline D } \right),\,0 < \alpha  < 2,\,s = \left[ \alpha  \right] + 1,\,k = \overline {0,s - 1} ;$

3)$\mathop {\lim }\limits_{y \to  + 0} \frac{{{\partial ^k}}}{{\partial {y^k}}}\left( {I_{0y}^{\left( {1 - \beta } \right)\left( {s - \alpha } \right)}u\left( {x,y} \right)} \right) = {\varphi _k}\left( x \right),\,\,k = \overline {0,s - 1} ,\,s = \left[ \alpha  \right] + 1.$\\
Given functions satisfy the constraints:
\[{\varphi _k}\left( x \right) \in C\left( R \right);\]
\[\left| {{\varphi _k}\left( x \right)} \right| < M\exp \left( {N{{\left| x \right|}^{\frac{{2n}}{{2n - \alpha }}}}} \right),\,\left| x \right| \to  + \infty ,\]
$$N < \sigma, 0<M-constant.$$
Taking Lemma 2 into account, we obtain the following result.

\textbf{Theorem 2.} The solution to the Cauchy problem has the form
\[u\left( {x,y} \right) = \sum\limits_{k = 0}^{s - 1} {\int\limits_{ - \infty }^{ + \infty } {{\varphi _k}\left( \xi  \right){\Gamma _{ - \frac{\alpha }{{2n}} - \left( {1 - \beta } \right)\left( {s - \alpha } \right) + k}}\left( {x - \xi ,y} \right)d\xi } },\eqno(19)\]
$$s = \left[ \alpha  \right] + 1,\,0 < \alpha  < 2.$$

Note that solutions of the Cauchy problem for equation (10) with fractional Riemann-Liouville and Caputo derivatives are written out through formula (19) for $ \beta = 0 $ and $ \beta = 1 $, respectively.

The case $ p = 2n + 1 $ is considered similarly.

\begin{center}
\textbf{References}
\end{center}
 1. R. Hilfer, Applications of Fractional Calculus in Physics.// World Scientific, Singapore (2000).\\
 2.  L u c h k o Yu.,  G o r e n f l o R. Scale-invariant solutions of a partial
differential equation of fractional order. // Fractional Calculus and Applied
Analysis.1998. Т.1, No 1. pp. 63--78. \\
 3. Irgashev B. Y. Construction of Singular Particular Solutions Expressed via Hypergeometric Functions for an Equation with Multiple Characteristics. //Differential Equations. 2020. Т. 56. No 3, pp. 315-323.\\
 4. M.-Ha. Kim, G.Chol-Ri, Chol O.H. Operational method for solving multi-term
fractional differential equations with the generalized fractional derivatives.// Fractional Calculus and Applied Analysis. 2014. vol.17. No 1, pp.79-95.\\
5. E.Karimov. Tricomi type boundary value problem with integral conjugation condition for a mixed type equation with Hilfer fractional operator.// arXiv:1904.00258.\\
6.  M. S. Salakhitdinov, E. T. Karimov. Direct and inverse source problems for two-term time-fractional diffusion equation with Hilfer derivative.// arXiv:1711.00352.\\
7. Pskhu A. Fundamental solutions and cauchy problems for an odd-order partial differential equation with. // Electronic Journal of Differential Equations, Vol. 2019 (2019), No. 21, pp. 1–13. \\
8. Karasheva L. L. Cauchy problem for high even order parabolic equation with time fractional derivative //Sibirskie Elektronnye Matematicheskie Izvestiya [Siberian Electronic Mathematical Reports].2018. Т. 15. pp. 696-706. [In Russian].\\
9.  Wright E. M. The generalized Bessel function of order greater than one // Quart.
J. Math., Oxford Ser., 1940. Т. 1, No 11. pp. 36–48.\\
\end{document}